\documentclass[a4paper, 11pt]{article}  
\usepackage{a4wide}

\usepackage{amssymb, amsmath, amsfonts, amssymb, mathtools}
\usepackage{hyperref}
\usepackage{xcolor}
\usepackage{cite}
\usepackage{caption}
\usepackage{subcaption}

\renewcommand{\leq}{\leqslant}
\renewcommand{\geq}{\geqslant}

\newcommand{\R}{\mathbb{R}}
\newcommand{\N}{\mathbb{N}}

\DeclareMathOperator{\Id}{Id}

\newcommand{\dd}{\textrm{d}}

\newcommand{\XX}{X}
\newcommand{\YY}{Y}
\newcommand{\HH}{H}

\newcommand{\LL}{\mathcal{L}}
\newcommand{\OO}{\mathcal{O}}

\newcommand{\xhat}{\hat{x}}
\newcommand{\yhat}{\hat{y}}
\newcommand{\bhat}{\hat{B}}
\newcommand{\xeps}{\tilde{x}}
\newcommand{\yeps}{\tilde{y}}
\newcommand{\beps}{\tilde{B}}

\newcommand{\zhat}{\hat{z}}

\DeclareMathOperator{\Tr}{Tr}

\renewcommand{\epsilon}{\varepsilon}
\newcommand{\eps}{\epsilon}
\newcommand{\gfrak}{\mathfrak{g}}

\newcommand{\startmodif}{}
\newcommand{\stopmodif}{}

\DeclareMathOperator{\dom}{\mathcal{D}}

\newtheorem{theorem}{Theorem}[section]
\newtheorem{assumption}[theorem]{Assumption}
\newtheorem{definition}[theorem]{Definition}
\newtheorem{problem}[theorem]{Problem}
\newtheorem{remark}[theorem]{Remark}
\newtheorem{example}[theorem]{Example}
\newtheorem{corollary}[theorem]{Corollary}

\usepackage{authblk}

\newcommand{\startmodiff}{}
\newcommand{\stopmodiff}{}                                 

\title{\LARGE \bf
Online estimation of Hilbert-Schmidt operators and application to kernel reconstruction of neural fields
}

\author[1]{Lucas Brivadis}
\author[1]{Antoine Chaillet}
\author[1]{Jean Auriol}
\affil[1]{Université Paris-Saclay, CNRS, CentraleSupélec, Laboratoire des Signaux et Systèmes, 91190, Gif-sur-Yvette, France. Emails:
        {\tt\small lucas.brivadis@centralesupelec.fr, antoine.chaillet@centralesupelec.fr, jean.auriol@centralesupelec.fr}}%

\begin{document}

\maketitle

\begin{abstract}

An adaptive observer is designed for online estimation of Hilbert-Schmidt operators
from online measurement of 
\startmodiff
part of
\stopmodiff
the state
for some class of nonlinear infinite-dimensional dynamical systems.
Convergence is ensured under detectability and persistency of excitation assumptions.
The class of systems considered is motivated by an application to kernel reconstruction of neural fields, commonly used to model spatiotemporal activity of neuronal populations. Numerical simulations confirm the relevance of the approach.

\end{abstract}

\section{Introduction}

The problem of online estimation of unknown parameters in dynamical systems from measured state variables is a major issue in many control systems.
It can be addressed by means of adaptive observers, that are observers estimating the unmeasured part of the state and the unknown parameters simultaneously.
The theory of adaptive observer design, well-known for linear finite-dimensional systems (see, e.g., \cite{sastry1990adaptive}), is still an active area of research when it comes to nonlinear \cite{BESANCON2000271, BESANCON201715416, https://doi.org/10.48550/arxiv.2112.05497} and/or infinite-dimensional \cite{DEMETRIOU19965346, DEMETRIOU2018220, 650677} systems. In this paper, we design an adaptive observer for a class of nonlinear infinite-dimensional systems that allows the reconstruction of unknown linear operators appearing in the dynamics. These operators are estimated in the Hilbert-Schmidt topology. Therefore, not only the state of the system is infinite-dimensional, but also the ``parameters'' (now, operators) to be estimated.

The specific class of systems we consider is motivated by an application to kernel reconstruction in neural fields.
The offline estimation of these kernels 
is now a classical issue in inverse problems for neuroscience (see
\cite{alswaihli2018kernel, potthast2009inverse} and references therein), that can be addressed for instance using a Tikhonov regularization. We instead rely on adaptive
observer strategies to
address the online estimation problem. The crucial additional constraint is that the reconstruction can only be based on past values of the measurements and estimates. The recent work \cite{https://doi.org/10.48550/arxiv.2111.02176} considers a similar problem but uses finite-dimensional conductance-based
models (which differ from the infinite-dimensional Wilson-
Cowan type equation considered here), and estimate finite-dimensional
parameters (while we reconstruct linear operators on infinite-
dimensional spaces).

\paragraph*{Notation} Given a Hilbert space $X$, we denote by $\langle\cdot,\cdot\rangle_X$ and $\|\cdot\|_X$ its corresponding scalar product and norm.
The identity operator over $X$ is denoted by $\Id_X$.
If $Y$ is a Hilbert space, we denote by $\LL(X, Y)$ the space of bounded linear operators from $X$ to $Y$ endowed with the operator norm 
\startmodiff
$\|\cdot\|_{\LL(X, Y)}$,
\stopmodiff
\startmodif
and we set $\LL(X) = \LL(X, X)$.
\stopmodif
For all $B\in\LL(X, Y)$, we denote by $B^*\in\LL(Y, X)$ its adjoint. If $B\in\LL(X)$, we denote by $\Tr(B)$ the trace of $B$ if it exists.
\startmodif
The operator $P\in\LL(X)$ is said to be self-adjoint positive-definite if $P=P^*$ and $\langle Px, x\rangle_X>0$ for all $x\in X\setminus\{0\}$.
\stopmodif
For any open interval $I\subset\R$,
\startmodiff
any $m\in\N$ and any $p\in[1, +\infty]$,
\stopmodiff
$L^p(I, X)$ and $W^{m, p}(I, X)$ stand for the usual Lebesgue and Sobolev spaces, endowed with their canonical norms.

\section{Problem statement}\label{sec:PS}

\subsection{Functional setting}\label{sec:functional}
Let $(X, \|\cdot\|_X)$ and $(Y, \|\cdot\|_Y)$ be two separable Hilbert spaces.
Consider an infinite-dimensional dynamical system of the following form:
\begin{equation}\label{eq:syst}
\left\{
\begin{aligned}
    &\dot x = A_1(x) + \psi(y) + u_1\\
    &\dot y = A_2(y) + B_1\phi_1(x) + B_2\phi_2(y) + u_2
\end{aligned}
\right.
\end{equation}
where $(x, y)$ is the state of the system lying in $X\times Y$, $u_1$ and $u_2$ are inputs respectively lying in $X$ and $Y$, and
$A_1:\dom(A_1)\to X$ and $A_2:\dom(A_2)\to Y$ are
singled-valued $m$-dissipative operators 
\startmodiff
(see \cite[Chapter 2]{miyadera1992nonlinear} for a definition),
\stopmodiff
respectively defined on dense subsets $\dom(A_1)\subset X$ and  $\dom(A_2)\subset Y$, such that $A_1(0) = 0$ and $A_2(0) = 0$. The linear operators $B_1\in\LL(X, Y)$ and $B_2\in\LL(Y)$ are bounded, and
$\psi:Y\to X$, $\phi_1:X\to X$ and $\phi_2:Y\to Y$ are Lipschitz continuous on any bounded set.
According to \cite[Chapter IV, Proposition 3.1]{showalter_monotone_2013},
$A_1$ and $A_2$ are generators of nonlinear strongly continuous contraction semigroups over $X$ and $Y$ respectively.


It follows from \cite[Chapter IV, Theorems 4.1 and 4.1A]{showalter_monotone_2013} that if $u_1$ and $u_2$ are absolutely continuous over $\R_+$, then for all $(x_0, y_0)\in \dom(A_1)\times\dom(A_2)$, there exists $t_{\max}\in (0,+\infty]$ such that \eqref{eq:syst} admits a unique strong solution $(x, y):[0,t_{\max})\to X\times Y$, i.e., such that $(x(0), y(0)) = (x_0, y_0)$, $(x, y)$ is absolutely continuous, satisfies \eqref{eq:syst} almost everywhere and lies in $\dom(A_1)\times\dom(A_2)$.
Moreover, if $(x, y)$ is bounded in $X\times Y$ over $[0,t_{\max})$, then $t_{\max} = +\infty$.

\subsection{Problem formulation}

In this paper, we consider the following online estimation problem.
\begin{problem}\label{problem}
From the knowledge of $A_1$, $A_2$, $\psi$, $\phi_1$, $\phi_2$
and the online measurement of $u_1$, $u_2$ and $y$, estimate online $x$ and the operators $B_1$ and $B_2$.
\end{problem}

In addition to the hypotheses made to ensure the well-posedness of the system,
we consider the following two main assumptions.

\begin{assumption}[Strong dissipativity]\label{ass:diss}
The nonlinear operator $A_1$ is strongly dissipative, that is, there exists a positive constant $\alpha$ such that for all $(x_1, x_2)\in \dom(A_1)^2$,
\begin{equation}
    \langle A_1(x_1) - A_1(x_2), x_1-x_2\rangle_X \leq -\alpha\|x_1-x_2\|^2_X.
\end{equation}
\end{assumption}
Since $A_1$ is supposed to be $m$-dissipative, we already have that $\langle A_1(x_1) - A_1(x_2), x_1-x_2\rangle_X \leq 0$, so that Assumption~\ref{ass:diss} is indeed a stronger dissipativity assumption.
Assumption~\ref{ass:diss} implies that any two solutions of $\dot x = A_1(x)$ are exponentially converging to one another in $X$ at exponential rate $\alpha$.
Since $y$ is supposed to be known online while $x$ is unknown, Assumption~\ref{ass:diss} can be interpreted as a detectability hypothesis: the unknown part of the state has a contracting dynamics.
In the estimation strategy, this allows to estimate $x$ online simply by simulating a particular trajectory of the $x$-subsystem (as all other solutions will eventually converge to it).
Note that,
\startmodiff
in some applications, 
\stopmodiff
the unmeasured state $x$ can also be assumed of null dimension
\startmodiff
(meaning full state
measurement)
\stopmodiff
, in which case the system \eqref{eq:syst} reduces to
\begin{equation*}
    \dot y = A_2(y) + B_2\phi_2(y) + u_2.
\end{equation*}
All the results of the paper are still valid in that easier case.

\startmodif

\begin{definition}
Let $(V_1, \|\cdot\|_{V_1})$ and $(V_2, \|\cdot\|_{V_2})$ be two separable Hilbert spaces. The linear bounded operator $B\in\LL(V_1, V_2)$ is said to be Hilbert-Schmidt if for any Hilbert basis $(e_k)_{k\in\N}$ of $(V_1, \|\cdot\|_{V_1})$,
\begin{align}\label{eq:hs}
    \|B\|^2_{\LL_2((V_1, \|\cdot\|_{V_1}), (V_2, \|\cdot\|_{V_2}))} :=\sum_{k\in\N} \|Be_k\|_{V_2}^2<+\infty.
\end{align}
\end{definition}
We denote by $\LL_2((V_1, \|\cdot\|_{V_1}), (V_2, \|\cdot\|_{V_2}))$ the Hilbert space of Hilbert-Schmidt operators from $V_1$ to $V_2$, endowed with the norm defined in \eqref{eq:hs}. When $V_1$ and $V_2$ are endowed with their norms $\|\cdot\|_{V_1}$ and $\|\cdot\|_{V_2}$ respectively, we simply write $\LL_2(V_1, V_2)$. We set $\LL_2(V_1) = \LL_2(V_1, V_1)$.
By definition of the trace operator $\Tr$ and of the adjoint $B^*$ of $B$, we have $\|B\|^2_{\LL_2(V_1, V_2)} = \Tr(B^*B) = \|B^*\|^2_{\LL_2(V_2, V_1)}$. 
Note that the topology induced by the Hilbert-Schmidt norm is finer than the one induced by the operator norm since $\|\cdot\|_{\LL(V_1, V_2)}\leq\|\cdot\|_{\LL_2(V_1, V_2)}$.

If $P\in\LL(V_1)$ is a self-adjoint positive-definite operator, then $\langle P\cdot, P\cdot\rangle_{V_1}$ defines a new scalar product on $V_1$, whose associated norm $\|P\cdot\|_{V_1}$ is weaker than or equivalent to $\|\cdot\|_{V_1}$.
Then, for all $B\in\LL_2(V_1, V_2)$,
$\|BP\|_{\LL_2(V_1, V_2)} = \|PB^*\|_{\LL_2(V_1, V_2)} = \|B^*\|_{\LL_2((V_1, \|\cdot\|_{V_1}), (V_2, \|P\cdot\|_{V_2}))}$. Thus $\|\cdot P\|_{\LL_2(V_1, V_2)}$ defines a norm on 
$\LL_2(V_1, V_2)$ that is weaker than or equivalent to $\|\cdot\|_{\LL_2(V_1, V_2)}$.

\begin{assumption}[Hilbert-Schmidt operators]\label{ass:hilbert}
The linear bounded operators $B_1$ and $B_2$ are in $\LL_2(X, Y)$ and $\LL_2(Y)$, respectively.
\end{assumption}

Then, Problem~\ref{problem} consists in finding $\xhat(t)$, $\hat B_1(t)$ and $\hat B_2(t)$ for all $t\geq0$
such that
$\|\xhat(t) - x(t)\|_{X}\to0$,
$\|\hat B_1(t) - B_1\|_{\LL_2(X, Y)}\to0$
and $\|\hat B_2(t) - B_2\|_{\LL_2(Y)}\to0$
as $t\to+\infty$ 
(when $X$ and $Y$ are endowed with some norms).
Such estimators must only depend at time $t$ on the knowledge of  $A_1$, $A_2$, $\psi$, $\phi_1$, $\phi_2$, $u_1(s)$, $u_2(s)$ and $y(s)$ for $s\in[0, t]$.
\stopmodif

\section{Kernel reconstruction of neural fields}
\label{sec:KR}

\subsection{Neural fields}

Problem~\ref{problem} is motivated by an application to kernel reconstruction of neural fields.
Neural fields are nonlinear integro-differential equations modeling the spatiotemporal evolution of the activity of neuronal populations. They are based on the seminal works \cite{wilson1973mathematical, amari1977dynamics} and surveys on their extensive use in mathematical neuroscience can be found in \cite{Bressloff2011,coombes2014neural}.
Given a compact set $\Omega\subset \mathbb R^q$ (where, typically, $q\in\{1,2,3\}$) representing the physical support of the population, the evolution of neuronal activity $z(t, r)\in\R^n$ at time $t\in\R_+$ and position $r\in\Omega$ is modeled as
\begin{align}\label{eq:wc}
    \tau(r)\frac{\partial z}{\partial t}(t, r) = -z(t, r) + u(t, r)+ \int_{r'\in\Omega}w(r, r')S(z(t, r'))\dd r'
    ,
\end{align}
where
\startmodiff
$n\in\N$ represents the number of considered neuronal
population types,
$\tau(r)$ is a positive diagonal matrix of size $n\times n$ continuous in $r$ representing the time decay constant of neuronal activity at position $r$,
\stopmodiff
$S:\R^n\to\R^n$ is a nonlinear activation function (typically, a sigmoid), $w(r, r')\in\R^{n\times n}$ defines a kernel describing the synaptic strength between location $r$ and $r'$ and $u(t, r)\in\R^n$ is an input.
We consider the problem of online reconstruction of the kernel $w$ from the measurement of the neuronal activity $z$.

\subsection{Application}\label{sec:neural}

Now we show how \eqref{eq:wc} fits into \eqref{eq:syst} and discuss the relevance of Assumptions \ref{ass:diss} and \ref{ass:hilbert} in this context.
In order to ensure well-posedness, we make the following usual assumptions on $S$ and $w$:
\begin{itemize}
    \item $S$ is bounded, differentiable and has bounded derivative;
    \item $w$ is square-integrable over $\Omega^2$.
\end{itemize}
These assumptions are standard in neural fields analysis. In particular, the boundedness of $S$ reflects the biological limitations of the maximal activity that can be reached by the population.

We assume that the neuronal population can be decomposed into $z(t, r)=(z_1(t, r), z_2(t, r))\in\R^{n-m}\times \R^m$ where $z_1$ corresponds to the unmeasured part of the state and $z_2$ to the measured part. Such a decomposition is natural when the two considered populations are physically separated, as it happens in the brain structures involved in Parkinson's disease \cite{DECH16}. It can also be relevant for imagery techniques that discriminate among neuron types within a given population. Accordingly, we define $\tau_i$, $S_i$ $w_{ij}$ and $u_i$ of suitable dimensions for each population $i, j\in\{1,2\}$ so that
\begin{align}\label{eq:wcij}
    \tau_i(r)\frac{\partial z_i}{\partial t}(t, r) = -z_i(t, r) + u_i(t, r)+ \sum_{j=1}^2\int_{r'\in\Omega}w_{ij}(r, r')S_j(z_j(t, r'))\dd r'
    .
\end{align}
Denote by $m$ the dimension of the measured activity $z_2(t,r)$. In order to fit \eqref{eq:wcij} in the form of \eqref{eq:syst}, set
$X = L^2(\Omega; \R^{n-m})$, $Y = L^2(\Omega; \R^{m})$,
\startmodiff
$x = z_1$, $y = z_2$,
$W_{ij}(z_j)=\textstyle\int_{r'\in\Omega}w_{ij}(\cdot, r')z_j(r')\dd r'$,
$A_1 = \tau_1^{-1}(-\Id_X + W_{11}S_1)$, $\dom(A_1) = X$,
$A_2 = -\tau_2^{-1}\Id_Y$, $\dom(A_2) = Y$,
$\psi = \tau_1^{-1}W_{12}S_2$,
$\phi_j = S_j$ and
$B_j = \tau_2^{-1}W_{2j}$.
\stopmodiff

Since $w$ is square-integrable, $B_1$ and $B_2$ are Hilbert-Schmidt integral operators with kernels $\tau_2^{-1}w_{21}$ and $\tau_2^{-1}w_{22}$, hence Assumption~\ref{ass:hilbert} is satisfied.
In order to satisfy the detectability Assumption~\ref{ass:diss}, we need to assume that $z_1$ has a strongly dissipative
\startmodiff
internal
\stopmodiff
dynamics, namely, that $A_1$ is strongly dissipative.
Remark that due to the structure of $A_1$, this is the case if
\begin{align}\label{eq-1}
    \ell_1\|W_{11}\|_{\LL(X)}<1
\end{align}
where $\ell_1$ is the Lipschitz constant of $S_1$. Indeed, it yields
\startmodiff
\begin{align*}
    \langle &A_1(x_1) - A_1(x_2), x_1-x_2\rangle_X\\
    &
    =- \|\tau_1^{-1/2}(x_1-x_2)\|_X^2
    + \langle W_{11}(S_1(x_1) - S_1(x_2)), \tau_1^{-1}(x_1-x_2)\rangle_X
    \\
    &\leq -\alpha\|x_1-x_2\|_X^2
\end{align*}
for $\alpha = \frac{1-\ell_1\|W_{11}\|_{\LL(X)}}{\max_\Omega \tau_1}$.
\stopmodiff
We stress that condition \eqref{eq-1} is commonly used in the stability analysis of neural fields \cite{Faugeras:2008wx} and ensures dissipativity even in the presence of axonal propagation delays \cite{DECHcdc17}.

We thus assume that each population is either measured online (taken into account in $z_2$) or unmeasured but 
\startmodiff
internally
strongly dissipative
\stopmodiff
and with known kernels (taken into account in $z_1$).
Problem~\ref{problem} is now equivalent to online reconstruction of $w_{21}$ and $w_{22}$ (in $L^2(\Omega^2)$) from the online measurement of $z_2$ and $u_i$ and the knowledge of $\tau_i$, $w_{1i}$, $S_i$ for all $i\in\{1,2\}$. Note that if the full state $z$ is measured (i.e. $m=n$), then no dissipative part $z_1$ of the system is required, hence the full kernel $w$ is to be estimated.

\section{Online estimation of Hilbert-Schmidt operators}
\label{sec:OE}

\subsection{Adaptive observer design}

In order to solve Problem~\ref{problem}, we propose to consider $B_1$ and $B_2$ as additional constant variables to system \eqref{eq:syst}, so that the resulting state space is the Hilbert space $H:=X\times Y\times \LL_2(X, Y)\times\LL_2(Y)$. Set also
$\dom:=\dom(A_1)\times \dom(A_2)\times \LL_2(X, Y)\times\LL_2(Y)\subset H$.
Inspired by the estimator proposed in \cite{BESANCON2000271} for finite-dimensional nonlinear systems, we consider the following observer over $H$:
\begin{equation}\label{eq:obs}
\left\{
\begin{aligned}
    &\dot \xhat = A_1(\xhat) + \psi(y) + u_1\\
    &\dot \yhat = A_2(y) + \bhat_1\phi_1(\xhat) + \bhat_2\phi_2(y) + u_2 - \beta(\yhat-y)\\
    &\dot{\bhat}_1 = - \gamma_1(\yhat-y)\phi_1(\xhat)^*\\
    &\dot{\bhat}_2 = -\gamma_2(\yhat-y)\phi_2(y)^*
\end{aligned}
\right.
\end{equation}
where $\beta$, $\gamma_1$, and $\gamma_2$ are positive constants, called observer gains, that need to be appropriately tuned to guarantee the convergence of the observer state to the real state.
Note that for any $v$ in $Y$ and any $w$ in $X$ (resp. in $Y$), $vw^*$ lies in $\LL_2(X, Y)$ (resp. in $\LL_2(Y)$) and $\|vw^*\|_{\LL_2(X, Y)} = \|v\|_Y\|w\|_X$ (resp. $\|vw^*\|_{\LL_2(Y)} = \|v\|_Y\|w\|_Y$).
Reasoning as in Section~\ref{sec:functional}, one can show that the cascade system \eqref{eq:syst}-\eqref{eq:obs} is well-posed.

\subsection{Main result}

Our main result, proved in Section \ref{sec-proof-theo}, relies on the notion of persistence of excitation.

\startmodif
\begin{definition}[Persistence of excitation]\label{ass:pe}
A continuous signal $g:\R_+\to V$ is persistently exciting over the Hilbert space $(V, \|\cdot\|_V)$ with respect to a self-adjoint positive-definite operator $P\in\LL(V)$ if there exists positive constants $T$ and $\kappa$ such that
\begin{equation}\label{eq:pe}
    \int_t^{t+T}g(\tau)g(\tau)^*\dd\tau \geq \kappa P^2,\quad \forall t\geq0.
\end{equation}
\end{definition}
\smallskip
\begin{remark}
If $V$ is finite dimensional, then Definition~\eqref{ass:pe} coincides with the usual notion of persistence of excitation since all norms on $V$ are equivalent.
However, if $V$ is infinite-dimensional, then there do not exist any persistently exciting signal with respect to the identity operator on $V$. (Actually, it is a characterization of the infinite-dimensionality of $V$). Indeed, if $P=\Id_V$, then \eqref{eq:pe} at $t=0$ together with the spectral theorem for compact operators implies that $\int_0^{T}g(\tau)g(\tau)^*\dd \tau$ is not a compact operator, which is in contradiction with the fact that the sequence of finite range operators $\sum_{j=0}^Ng(\frac{jT}{N})g(\frac{jT}{N})^*$ converges to it in $\LL(V)$ as $N$ goes to infinity.
This is the reason for which we introduce this new persistency of excitation condition which is feasible even if $V$ is infinite-dimensional.
\end{remark}

\begin{example}
Let $V=l^2(\N, \R)$ be the space of square summable real sequences. The signal $g:\R_+\to V$ defined by $g(\tau) = (\frac{\sin(k\tau)}{k^2})_{k\in\N}$ is \startmodiff
persistently exciting
\stopmodiff
with respect to $P\in\LL(V)$ defined by $P(x_k)_{k\in\N} = (\frac{x_k}{k^2})_{k\in\N}$ with constant $T=2\pi$ and $\kappa = \pi$ since $\int_0^{2\pi}\sin^2(k\tau)\dd\tau = \pi$.
\end{example}

\stopmodif

We now provide sufficient conditions for the convergence of the observer \eqref{eq:obs} to the state of system \eqref{eq:syst}, thus solving Problem \ref{problem}.

\begin{theorem}[Observer convergence]\label{th:main}
Suppose that Assumptions~\ref{ass:diss} and~\ref{ass:hilbert} are satisfied.
Assume moreover that the functions~$\phi_1$ and $\phi_2$ are bounded and that
$\phi_1$ is globally Lipschitz continuous with constant $\ell_1$.
Pick the observer gains $\beta,\gamma_1,\gamma_2$ such that  $\gamma_1,\gamma_2>0$ and
\begin{equation}\label{eq:ab}
4\alpha\beta>\ell^2_1\|B_1\|^2_{\LL(X, Y)}.
\end{equation}
Then, for any absolutely continuous $u_1$ and $u_2$ and any solution of \eqref{eq:syst} defined over $\R_+$,
any solution of \eqref{eq:obs} satisfies
\begin{align*}
    \lim_{t\to+\infty} \|\xhat(t) - x(t)\|_X=0, \quad \lim_{t\to+\infty}\|\yhat(t) - y(t)\|_Y=0,
\end{align*}
and
$\|\hat B_1 - B_1\|_{\LL_2(X, Y)}$
and
$\|\hat B_2 - B_2\|_{\LL_2(Y)}$
remain bounded.

\startmodif
Moreover, if $P_X\in\LL(X)$ and $P_Y\in\LL(Y)$ are self-adjoint positive-definite operators such that $t\mapsto(\phi_1(x(t)), \phi_2(y(t)))$ is persistently exciting over $X\times Y$ with respect to $P=\begin{pmatrix}P_X&0\\0&P_Y\end{pmatrix}\in\LL(X\times Y)$,
and if $\phi_1$ and $\phi_2$ are differentiable with bounded derivatives,
then
\begin{align*}
    &\lim_{t\to+\infty} \|(\hat B_1(t) - B_1)P_X\|_{\LL_2(X, Y)}=0,\\
    &\lim_{t\to+\infty} \|(\hat B_2(t) - B_2)P_Y\|_{\LL_2(Y)}=0.
\end{align*}
\stopmodif
\end{theorem}
It is worth noting that the observer gains $\gamma_1$ and $\gamma_2$ play no qualitative role in the observer convergence. Also,~$\beta$ can always be picked sufficiently large to fulfill \eqref{eq:ab}. The main requirement therefore lies in the persistence of excitation requirement, which is a common hypothesis to ensure convergence of adaptive observers (see for instance \cite{BESANCON2000271, FARZA20092292, sastry1990adaptive} in the finite-dimensional context and \cite{DEMETRIOU19965346, 761927} in the infinite-dimensional case).
Roughly speaking, it states that parameters to be estimated are sufficiently ``excited'' by the system dynamics. However, this assumption is difficult to check in practice since it depends on the trajectories of the system itself.
In Section~\ref{sec:num}, we choose in numerical simulations a persistently exciting input $(u_1, u_2)$ in order to generate persistence of excitation in the signal $(\phi_1(x), \phi_2(y))$. This strategy seems to be numerically efficient, but the theoretical analysis of the link between the persistence of excitation of $(u_1, u_2)$ and $(\phi_1(x), \phi_2(y))$ remains an open question, not only in the present work but also for general classes of adaptive observers.

\startmodiff
\subsection{Application to neural fields}
\stopmodiff

As developed in Section~\ref{sec:neural}, Theorem~\ref{th:main} directly applies to the neural fields context. 
With the notations of Section~\ref{sec:neural}, the adaptive observer takes the form
\begin{equation}\label{eq:obswc}
\left\{
\begin{aligned}
    &\tau_1\dot \zhat_1 = -\zhat_1 + W_{11}S_1(\zhat_1) + W_{12}S_2(z_2) + u_1\\
    &
    \tau_2 \dot \zhat_2 = -z_2 +
    \hat W_{21} S_1(\zhat_1) + \hat W_{22} S_2(z_2)
    + u_2
    - \startmodiff
    \tau_2\beta
    \stopmodiff
    (\zhat_2-z_2)
    \\
    &\dot{\hat{W}}_{21} = - \gamma_1(\zhat_2-z_2)S_1(\zhat_1)^*\\
    &\dot{\hat{W}}_{22} = -\gamma_2(\zhat_2-z_2)S_2(z_2)^*.
\end{aligned}
\right.
\end{equation}
Then, we have the following, which immediately follows from Theorem \ref{th:main} for this particular system.

\begin{corollary}[Neural fields estimation]\label{cor:main}
Suppose that $S_1$ (resp. $S_2$) is bounded, differentiable, and that its derivative is bounded by some $\ell_1>0$ (resp. $\ell_2>0$), and that $w$ is square-integrable over $\Omega^2$.
Assuming that \eqref{eq-1} is satisfied, pick the observer gains in such a way that $\gamma_1,\gamma_2>0$ and
\begin{equation}\label{eq-2}
\startmodiff
4\frac{1-\ell_1\|W_{11}\|_{\LL(X)}}{\max_\Omega \tau_1}\beta
>\ell^2_1\|\tau_2^{-1}W_{21}\|^2_{\LL(X, Y)}.
\stopmodiff
\end{equation}
Consider any solution of \eqref{eq:wcij} defined over $\R_+$ for some absolutely continuous inputs $u_1$ and $u_2$. Then any solution of  \eqref{eq:obswc} satisfies $
\lim_{t\to+\infty}\|\zhat_1(t) - z_1(t)\|_X=0,~ \lim_{t\to+\infty} \|\zhat_2(t) - z_2(t)\|_Y=0,$
and
$\|\hat W_{21} -  W_{21}\|_{\LL_2(X, Y)}$
and
$\|\hat W_{22} -  W_{22}\|_{\LL_2(X, Y)}$
remain bounded.

Moreover,
then
\startmodiff
if $P_X\in\LL(X)$ and $P_Y\in\LL(Y)$ are self-adjoint positive-definite operators such that $t\mapsto(S_1(z_1(t)), S_2(z_2(t)))$ is persistently exciting over $X\times Y$ with respect to $P=\begin{pmatrix}P_X&0\\0&P_Y\end{pmatrix}\in\LL(X\times Y)$,
then
\begin{align*}
    &\lim_{t\to+\infty} \|(\hat W_{21}(t) -  W_{21})P_X\|_{\LL_2(X, Y)}=0,
    \\
    &\lim_{t\to+\infty} \|(\hat W_{22}(t) -  W_{22})P_Y\|_{\LL_2(Y)}=0,
\end{align*}
\stopmodiff
\end{corollary}

Here again, provided that condition \eqref{eq-1} holds, $\beta$ can always be picked large enough to fulfill \eqref{eq-2}.

\subsection{Proof of Theorem~\ref{th:main}}\label{sec-proof-theo}

\startmodiff
Consider a solution $(x, y)$ of \eqref{eq:syst} 
\stopmodiff
and the corresponding solution $(\xhat, \yhat, \bhat_1, \bhat_2)$ of \eqref{eq:obs}.
The estimation error $(\xeps, \yeps, \beps_1, \beps_2) := (\xhat, \yhat, \bhat_1, \bhat_2) - (x, y, B_1, B_2)$ is ruled by:
\begin{equation}\label{eq:eps}
\left\{
\begin{aligned}
    &\dot \xeps = A_1(\xhat) - A_1(x)\\
    &\dot \yeps = \bhat_1\phi_1(\xhat)-B_1\phi_1(x) + \beps_2\phi_2(y) - \beta\yeps\\
    &\dot{\beps}_1 = - \gamma_1\yeps\phi_1(\xhat)^*\\
    &\dot{\beps}_2 = -\gamma_2\yeps\phi_2(y)^*.
\end{aligned}
\right.
\end{equation}

\subsubsection{Proof that $(\xeps, \yeps)\to0$}

We endow $H$ with the squared norm
$\|\cdot\|_H^2 = \|\cdot\|_X^2 + \|\cdot\|_Y^2 + \frac{1}{\gamma_1}\|\cdot\|_{\LL_2(\XX, \YY)}^2 + \frac{1}{\gamma_2}\|\cdot\|_{\LL_2(\YY, \YY)}^2$, which is equivalent to the squared norm induced by the Cartesian product $H=X\times Y\times \LL_2(X, Y)\times\LL_2(Y)$.
Given any initial state, denote by $t_{\max}\in(0,+\infty]$ the maximal time of existence of $(\xeps, \yeps, \beps_1, \beps_2)$. Using
$\bhat_1\phi_1(\xhat)-B_1\phi_1(x)=\beps_1\phi_1(\xhat)+B_1(\phi_1(\xhat)-\phi_1(x))$, we have almost everywhere on $[0,t_{\max})$
\begin{align*}
    \frac{\dd}{\dd t} \|(\xeps, \yeps, \beps_1, \beps_2)\|_\HH^2
    &
    = \langle A_1(\xhat) - A_1(x), \xeps\rangle_\XX
    - \beta \|\yeps\|_\YY^2
    + \langle \beps_2\phi_2(y), \yeps\rangle_\YY
    \\
    &\quad+ \langle \beps_1\phi_1(\xhat), \yeps\rangle_\YY
    + \langle B_1(\phi_1(\xhat)-\phi_1(x)), \yeps\rangle_\YY
    \\
    &\quad
    - \langle \yeps\phi_1(\xhat)^*, \beps_1\rangle_{\LL_2(\XX, \YY)}
    - \langle \yeps\phi_1(y)^*, \beps_2\rangle_{\LL_2(\XX, \YY)}.
\end{align*}
By Assumption~\ref{ass:diss}, $\langle A_1(\xhat) - A_1(x), \xeps\rangle_\XX \leq -\alpha \|\xeps\|_\XX^2$.
By definition of the Hilbert-Schmidt scalar product,
\begin{align*}
\langle \yeps\phi_1(\xhat)^*, \beps_1\rangle_{\LL_2(\XX, \YY)}
&= \Tr(\phi_1(\xhat)\yeps^*\beps_1)
= \Tr(\yeps^*\beps_1\phi_1(\xhat))\\
&= \langle \beps_1\phi_1(\xhat), \yeps\rangle_\YY
\end{align*}
and, similarly,
$\langle \yeps\phi_2(y)^*, \beps_2\rangle_{\LL_2(\XX, \YY)}
=\langle \beps_2\phi_2(y), \yeps\rangle_\YY.$
Hence $
    \frac{\dd}{\dd t} \|(\xeps, \yeps, \beps_1, \beps_2)\|_\HH^2
    \leq -\alpha \|\xeps\|_\XX^2
    - \beta \|\yeps\|_\YY^2\\
    + \langle B_1(\phi_1(\xhat)-\phi_1(x)), \yeps\rangle_\XX.$
By Cauchy-Schwartz inequality, for all $\eps>0$,
\begin{align*}
\langle B_1(\phi_1(\xhat)-\phi_1(x)), \yeps\rangle_\XX
\leq \ell_1 \|B_1\|_{\LL(\XX, \YY)} \left(\frac{\eps}{2}\|\xeps\|_\XX^2 + \frac{1}{2\eps}\|\yeps\|_\XX^2\right)
\end{align*}
where $\ell_1$ is the Lipschitz constant of $\phi_1$.
Pick $\eps = \frac{\alpha}{\ell_1 \|B_1\|_{\LL(\XX, \YY)}}+\frac{\ell_1 \|B_1\|_{\LL(\XX, \YY)}}{4\beta}>0$.
Using condition~\eqref{eq:ab},
we get that
$\mu_1 := \alpha - \ell_1 \|B_1\|_{\LL(\XX, \YY)}\eps/2>0$
and
$\mu_2 := \beta - \ell_1 \|B_1\|_{\LL(\XX, \YY)}/2\eps>0$.
Then
\begin{align*}
    \frac{\dd}{\dd t} \|(\xeps, \yeps, \beps_1, \beps_2)\|_\HH^2
    \leq -\mu_1 \|\xeps\|_\XX^2
    - \mu_2 \|\yeps\|_\YY^2.
\end{align*}
Thus $(\xeps, \yeps, \beps_1, \beps_2)$ remains bounded. Hence according to 
\startmodiff
\cite[Theorem 4.10]{miyadera1992nonlinear},
\stopmodiff
we obtain
$t_{\max}=+\infty$ i.e. the state $(\xhat, \yhat, \bhat_1, \bhat_2)$ is defined over $\R_+$.
Moreover, we have
$
\frac{\dd}{\dd t}\|\xeps\|_\XX^2 \leq 0
$
and
\begin{align*}
    \frac{\dd}{\dd t}\|\yeps\|_\YY^2
    \leq
    - \beta \|\yeps\|_\YY^2
    + \langle \beps_2\phi_2(y), \yeps\rangle_\YY
    +\langle\bhat_1\phi_1(\xhat), \yeps\rangle_\YY
    -\langle B_1\phi_1(x), \yeps\rangle_\YY
\end{align*}
which is bounded since $(\xeps, \yeps, \beps_1, \beps_2)$ is bounded, $B_1$ is constant, and $\phi_1$ and $\phi_2$ are bounded.
Hence, according to Barbalat's lemma,
$(\xeps(t), \yeps(t))\to0$ as $t\to+\infty$.

\subsubsection{Proof that $(\beps_1, \beps_2)\to0$}
Now assume that $t\mapsto(\phi_1(x(t)), \phi_2(y(t)))$ is persistently exciting over $X\times Y$
\startmodif
with respect to $P$,
\stopmodif
and that $\phi_1$ and $\phi_2$ are differentiable with bounded derivatives.
Note that $\bhat_1\phi_1(\xhat)-B_1\phi_1(x) = \bhat_1(\phi_1(\xhat)-\phi_1(x)) + \beps_1\phi_1(x)$.
Hence the error dynamics \eqref{eq:eps} can be written as
\begin{equation}
\left\{
\begin{aligned}
    &\dot \yeps(t) = \beps_1(t)g_1(t) + \beps_2(t)g_2(t) + f_0(t)\\
    &\dot{\beps}_1(t) = f_1(t)\\
    &\dot{\beps}_2(t) = f_2(t),
\end{aligned}
\right.
\end{equation}
where
$f_0(t) := \bhat_1(\phi_1(\xhat(t))-\phi_1(x(t))) - \beta\yeps$,
$f_1(t) := - \gamma_1\yeps(t)\phi_1(\xhat(t))^*$,
$f_2(t) := - \gamma_2\yeps(t)\phi_2(y(t))^*$,
$g_1(t) := \phi_1(x(t))$ and $g_2(t) := \phi_2(y(t))$ for all $t\geq0$. Since $(\xeps(t), \yeps(t))\to0$ as $t\to+\infty$, $\bhat_1$ is bounded, $\phi_1$ is globally Lipschitz and $\phi_1$ and $\phi_2$ are bounded, we get that $f_i(t)$ tends toward 0 as $t$ goes to $+\infty$ for all $i\in\{0,1,2\}$.
Set
$g := (g_1, g_2):\R_+\to X\times Y$ and
$f_{1,2},\, \beps : \R_+\to \LL_2(X\times Y, Y)$,
defined by
$f_{1,2}(t)(\zeta, \xi) = f_1(t)\zeta + f_2(t)\xi$.
Set
$\beps(t)(\zeta, \xi) = \beps_1(t)\zeta + \beps_2(t)\xi$ for all $(\zeta, \xi)\in X\times Y$ and all $t\geq0$,
so that $\dot \yeps(t) = \beps(t) g(t) + f_0(t)$.
\startmodif
Remark that
$
\|\beps(t)P\|_{\LL_2(X\times Y, Y)}^2
=
\|\beps_1(t)P_X\|_{\LL_2(X, Y)}^2
+
\|\beps_2(t)P_Y\|_{\LL_2(Y)}^2
$, so that it remains to show that
$\|\beps(t)P\|_{\LL_2(X\times Y, Y)}\to0$ as $t\to+\infty$ to conclude.
\stopmodif

Applying twice Duhamel's formula, we have for all $t, \tau\geq0$: $\yeps(t+\tau)=
    \yeps(t)
    + \beps(t)\int_0^{\tau}g(t+s)\dd s
    \\
    + \int_0^{\tau}\int_0^sf_{1,2}(t+\sigma)g(t+s)\dd\sigma \dd s
    + \int_0^{\tau}f_0(t+s)\dd s.$
Define $\OO(t, T):=\int_0^{T}\|\yeps(t+\tau)\|_Y^2\dd \tau$ for any $T>0$ and $t\geq 0$.
Since $\yeps(t)\to0$, $\OO(t, T)\to0$ as $t\to+\infty$.
Moreover,
\begin{align*}
    \OO(t, T)
    &=
    \int_0^{T}\Big\|\yeps(t)+ \int_0^{\tau}\int_0^sf_{1,2}(t+\sigma)g(t+s)\dd\sigma \dd s\\
    &\quad
    + \int_0^{\tau}f_0(t+s)\dd s\Big\|_Y^2\dd \tau
    +
    \int_0^{T}\Big\|\beps(t)\int_0^{\tau}g(t+s)\dd s\Big\|_Y^2\dd \tau\\&\quad
    +\int_0^{T}\Big\langle\yeps(t)+ \int_0^{\tau}\int_0^sf_{1,2}(t+\sigma)g(t+s)\dd\sigma \dd s\\&\quad
    + \int_0^{\tau}f_0(t+s)\dd s
    ,\beps(t)\int_0^{\tau}g(t+s)\dd s\Big\rangle_Y\dd \tau.
\end{align*}
Since $\yeps(t)$ and $f_i(t)$ tends toward 0 as $t$ goes to $+\infty$ for all $i\in\{0,1,2\}$ and $g$ and $\beps$ are bounded, we get that
\begin{equation}\label{eq:conv1}
    \lim_{t\to+\infty}\int_0^{T}\Big\|\beps(t)\int_0^{\tau}g(t+s)\dd s\Big\|^2_Y\dd \tau=0.
\end{equation}
For all $t\geq 0$, define $h(t, \tau) = \beps(t)\int_0^{\tau}g(t+s)\dd s$.
By \eqref{eq:conv1}, $\|h(t, \cdot)\|_{L^2((0, T); Y)}\to0$ as $t\to+\infty$.
Note that $\frac{\partial h}{\partial \tau}(t, \tau) = \beps(t)g(t+\tau)$ hence $h(t, \cdot)\in W^{1, 2}((0, T); Y)$ since $g$ is bounded.
Moreover, $\dot \yeps$ is bounded since $\beps_i$ and $g_i$ are bounded for $i\in\{1,2\}$
and
$\dot\xeps = A_1(\xhat) - A_1(x)$ is bounded since $A_1$ is $m$-dissipative (see \cite[Corollary 3.7 and Theorem 4.20]{miyadera1992nonlinear}).
Hence,
if $\phi_1$ and $\phi_2$ are differentiable with bounded derivatives,
then so is $g$. Therefore, for all $t\geq 0$, $h(t, \cdot)\in W^{2, 2}((0, T); Y)$
and $\|h(t, \cdot)\|_{W^{2, 2}((0, T); Y)}\leq C_1$ for some positive constant $C_1$ independent of $t$.
According to the interpolation inequality (see, e.g., \cite[Section II.2.1]{temam1986infinite}),
\begin{equation*}
    \|h(t, \cdot)\|_{W^{1, 2}((0, T); Y)}^2\leq C_2\|h(t, \cdot)\|_{L^2((0, T); Y)}
\end{equation*}
for some positive constant $C_2$ independent of $t$.
Thus $\|\frac{\partial h}{\partial \tau}(t, \tau)\|_{L^2((0, T); Y)}\to0$, meaning that
\begin{equation}\label{eq:conv2}
    \lim_{t\to+\infty}\int_0^{T}\Big\|\beps(t)g(t+\tau)\dd s\Big\|^2_Y=0.
\end{equation}
Now, let $(e_k)_{k\in\N}$ be a Hilbert basis of $\YY$.
\startmodif
Then $\|\beps(t)P\|_{\LL_2(X\times Y, Y)}^2=\sum_{k\in\N} \|P\beps(t)^*e_k\|_{X\times Y}^2.$
Since $g$ is persistently exciting,
we have, for some $T,\kappa>0$,
\begin{align*}
    \int_0^{T}|\langle g(t+\tau), v\rangle_{X\times Y}|^2\dd\tau \geq \kappa\|Pv\|_{X\times Y}^2,\quad \forall t\geq0,
\end{align*}
for all $v\in \XX\times\YY$.
\stopmodif
Then,
\begin{align*} 
    \kappa \|\beps(t)P\|_{\LL_2(X\times Y, Y)}^2
    &\leq \sum_{k\in\N} \int_0^{T}|\langle g(t+\tau), \beps(t)^*e_k\rangle_{X\times Y}|^2\dd\tau\\
    &=\int_0^{T}\left\| \beps(t)g(t+\tau) \right\|_{Y}^2\dd\tau.
\end{align*}
Thus, by \eqref{eq:conv2},
$\|\beps(t)P\|_{\LL_2(X\times Y, Y)}^2\to 0$ as $t\to+\infty$, which concludes the proof.
\hfill$\blacksquare$

\startmodif
\section{Numerical simulation of kernel reconstruction of neural fields}\label{sec:num}
\stopmodif

We provide a numerical simulation of the adaptive observer 
\startmodiff
\eqref{eq:obswc}
\stopmodiff
in the case of a two-dimensional neural field (namely, $n=2$ and $m=1$ in Section~\ref{sec:neural}) over the unit circle $\Omega = \mathbb{S}^1$.
\startmodif
We set parameters of system \eqref{eq:wcij} and observer \eqref{eq:obswc} as in Table~\ref{table1}, so that all assumptions of Corollary~\ref{cor:main} are satisfied.
\stopmodif
Initial conditions are given by
\startmodiff
$z_1(0, r) = z_2(0, r) = 1$, $\zhat_1(0, r) = \zhat_2(0, r) = 0$ for all $r\in\Omega$
and $\hat W_{21}(0) = \hat W_{22}(0)=0$.
\stopmodiff
Kernels are given by Gaussian functions depending on the distance between $r$ and $r'$, as it is frequently assumed in practice (see \cite{DECH16}):
$w_{ij}(r, r') = \omega_{ij}\gfrak(r, r')/\|\gfrak\|_{L^2(\Omega^2; \R)}$, $\gfrak(r, r') = \exp(-\sigma|r-r'|^2)$ for constant parameters $\sigma$ and $\omega_{ij}$ given in Table~\ref{table1}.
The inputs $u_i$ are chosen as spatiotemporal periodic signals with irrational frequency ratio, i.e., $u_i(t, r) = 10^3\sin(\lambda_itr)$ with $\lambda_1/\lambda_2$ irrational. This choice is made to ensure persistency of excitation of the input $(u_1, u_2)$, which in practice seems to be sufficient to ensure persistency of excitation of $(S_1(z_1), S_2(z_2))$. Note that for $u_1=u_2=0$, the persistency of excitation assumption seems to be not guaranteed, hence the observer does not converge
\startmodiff
(the plot is not reported here).
\stopmodiff
However, in practice, such a persistent input is likely to occur due to exogenous signals coming from other unmodeled neuronal populations.
\startmodiff
Simulations code
\stopmodiff
can be found in repository \cite{git}. The system is spatially discretized over $\Omega$ with a constant space step $\Delta r = 1/20$, and the resulting ordinary differential equation is solved with an explicit Runge-Kutta $(4,5)$ method.

\begin{table}[ht!]
    \centering
    \vspace{0.25cm}
    \begin{tabular}{|c|c|c|c|}
        \hline
        $S_i(z) = \tanh(z)$&
        $\tau_i = 1$ &
        $\lambda_1 = 1$ &
        $\lambda_2 = \sqrt{2}$\\
        \hline
        $\omega_{11} = 0.1$&
        $\omega_{12} = 2$&
        $\omega_{21} = -2$&
        $\omega_{22} = 2$\\
        \hline
        $\beta = 1$&
        $\gamma_1 = 100$ & $\gamma_2 = 100$&
        $\sigma = 60$\\
        \hline
    \end{tabular}
    \caption{System and observer parameters for the numerical simulation of Figures~\ref{fig:Nx20c1000} and~\ref{fig:film}}
    \label{table1}
\end{table}

\begin{figure}[ht!]
    \centering
    \includegraphics[width=0.76\linewidth]{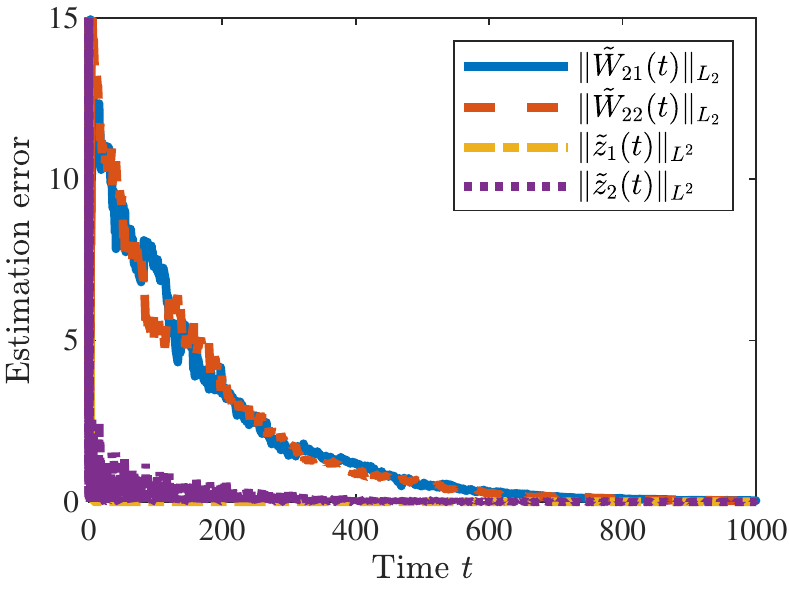}
    \caption{
    \startmodiff
    Evolution of the estimation errors 
$\|\hat W_{2i}(t) - W_{2i}\|_{\LL_2(X, Y)}$ and
$\|\hat z_i(t) - z_i(t)\|_X$
for $i\in\{1,2\}$.
\stopmodiff
}
\label{fig:Nx20c1000}
\end{figure}

\begin{figure}[ht!]
    \centering
    \begin{subfigure}[b]{0.49\linewidth}
         \centering
         \includegraphics[width=\linewidth]{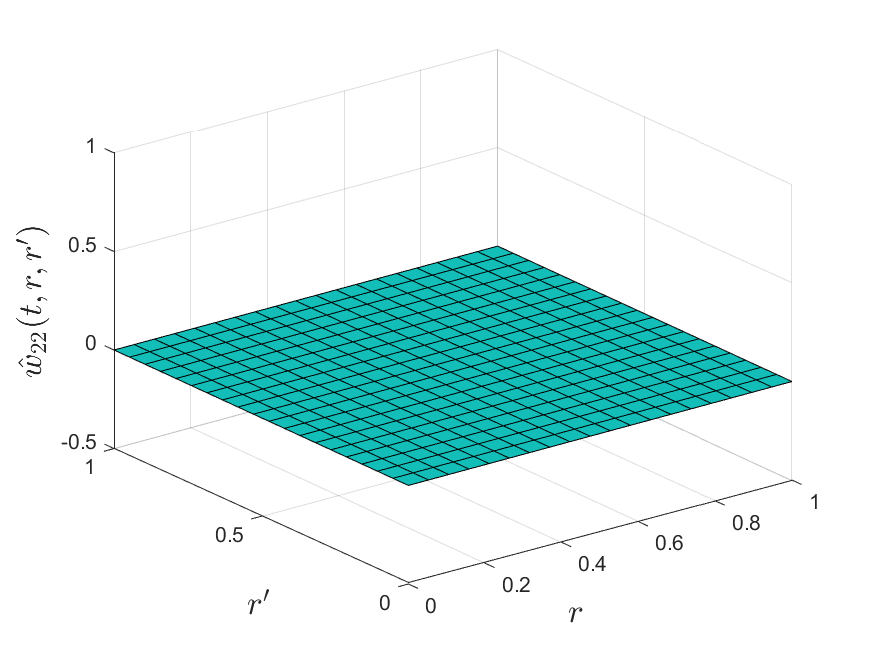}
         \caption{$t = 0$}
     \end{subfigure}
     \hfill
     \begin{subfigure}[b]{0.49\linewidth}
         \centering
         \includegraphics[width=\linewidth]{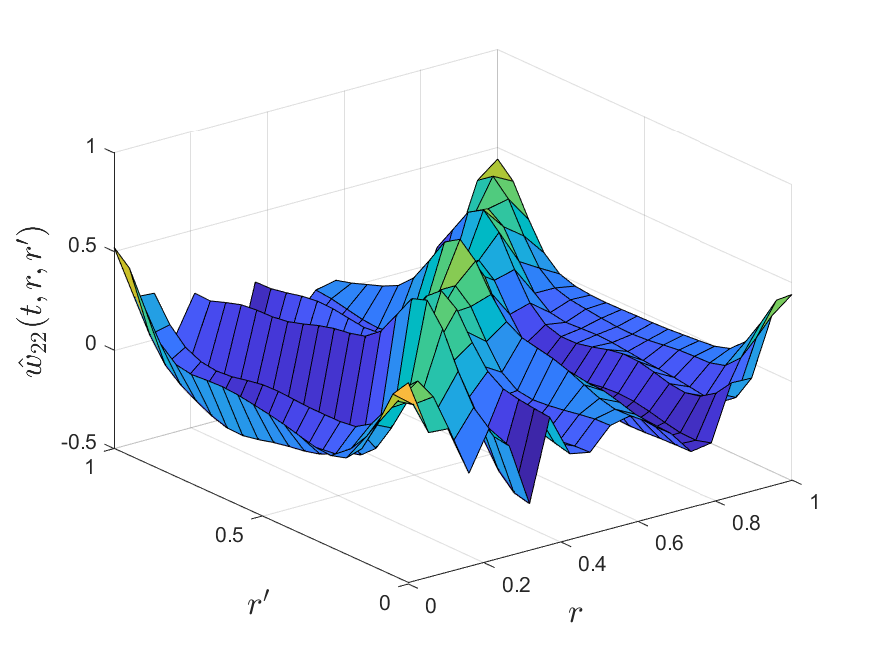}
         \caption{$t = 250$}
     \end{subfigure}
     \\
     \begin{subfigure}[b]{0.49\linewidth}
         \centering
         \includegraphics[width=\linewidth]{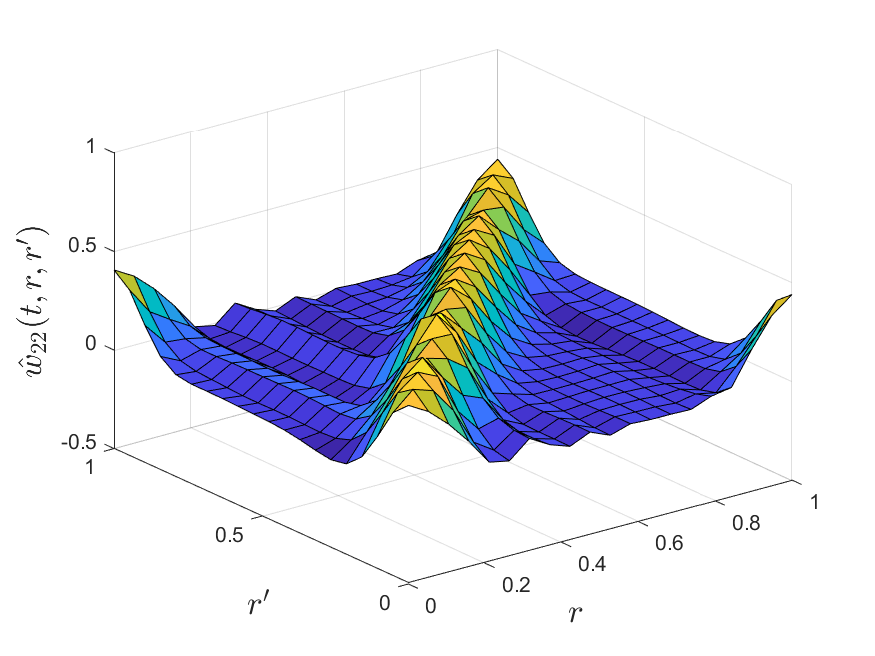}
         \caption{$t = 500$}
     \end{subfigure}
     \hfill
     \begin{subfigure}[b]{0.49\linewidth}
         \centering
         \includegraphics[width=\linewidth]{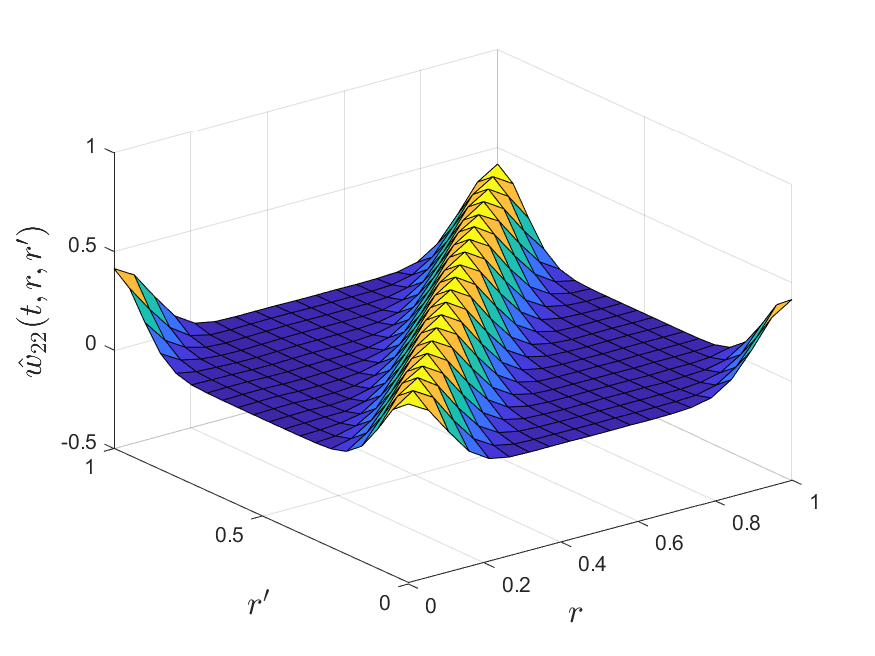}
         \caption{$t = 1000$}
     \end{subfigure}
    \caption{
    \startmodiff
    Evolution of the kernel $\hat w_{22}(t, r, r')$.
    \stopmodiff
    }
    \label{fig:film}
\end{figure}

In Figure~\ref{fig:Nx20c1000}, the convergence of the observer, that is proved in Corollary~\ref{cor:main}, is numerically verified. In Figure~\ref{fig:film}, we illustrate some iterations of the reconstructed kernel $\hat w_{22}(t)$ of
\startmodiff
$\hat W_{22}(t)$,
\stopmodiff
which converges to the kernel $w_{22}$ of $W_{22}$.

\section{Conclusion}

In this paper, we have shown that an observer can be designed to estimate online linear operators arising in some nonlinear infinite-dimensional dynamical systems from the measurement
\startmodiff
of part of the state variables, provided that the other variables have a strongly dissipative internal dynamics.
\stopmodiff
This estimation problem is motivated by an application to kernel reconstruction for neural fields equations.
The main assumption is the persistence of excitation of the system along its trajectories.
Our simulations suggest that this requirement can be ensured using appropriate exogenous inputs. In future works, we wish to investigate this hypothesis by either designing inputs ensuring \startmodiff
persistency of excitation,
\stopmodiff
or designing observers that do not rely on this assumption (see, e.g., \cite{8786148}).
Moreover, the use of this estimator in closed-loop to stabilize the systems by means of dynamic output feedback will be investigated.
Finally, delayed neural fields in the form of \cite{DECH16} could be considered, as they do not fit into the functional setting of the present paper although capturing meaningful biological processes such as non-instantaneous axonal propagation.

\bibliographystyle{abbrv}
\bibliography{references.bib}

\begin{thebibliography}{10}

\bibitem{alswaihli2018kernel}
J.~Alswaihli, R.~Potthast, I.~Bojak, D.~Saddy, and A.~Hutt.
\newblock Kernel reconstruction for delayed neural field equations.
\newblock {\em The Journal of Mathematical Neuroscience}, 8(1):1--24, 2018.

\bibitem{amari1977dynamics}
S.~Amari.
\newblock Dynamics of pattern formation in lateral-inhibition type neural
  fields.
\newblock {\em Biological cybernetics}, 27(2):77--87, 1977.

\bibitem{BESANCON2000271}
G.~Besançon.
\newblock Remarks on nonlinear adaptive observer design.
\newblock {\em Systems \& Control Letters}, 41(4):271--280, 2000.

\bibitem{BESANCON201715416}
G.~Besançon and A.~Ţiclea.
\newblock On adaptive observers for systems with state and parameter
  nonlinearities.
\newblock {\em IFAC-PapersOnLine}, 50(1):15416--15421, 2017.
\newblock 20th IFAC World Congress.

\bibitem{Bressloff2011}
P.~Bressloff.
\newblock {Spatiotemporal dynamics of continuum neural fields}.
\newblock {\em Journal of Physics A: Mathematical and Theoretical}, 45(3),
  2012.

\bibitem{git}
L.~Brivadis.
\newblock Kernelestimation project.
\newblock https://github.com/brivadis/KernelEstimation, 2021.

\bibitem{https://doi.org/10.48550/arxiv.2111.02176}
T.~B. Burghi and R.~Sepulchre.
\newblock Online estimation of biophysical neural networks, 2021.

\bibitem{DECH16}
A.~Chaillet, G.~Detorakis, S.~Palfi, and S.~Senova.
\newblock {Robust stabilization of delayed neural fields with partial
  measurement and actuation}.
\newblock {\em Automatica}, 83:262--274, Sep. 2017.

\bibitem{coombes2014neural}
S.~Coombes, P.~beim Graben, R.~Potthast, and J.~Wright.
\newblock {\em Neural Fields: Theory and Applications}.
\newblock Springer, 2014.

\bibitem{650677}
R.~Curtain, M.~Demetriou, and K.~Ito.
\newblock Adaptive observers for structurally perturbed infinite dimensional
  systems.
\newblock In {\em Proceedings of the 36th IEEE Conference on Decision and
  Control}, volume~1, pages 509--514 vol.1, 1997.

\bibitem{761927}
R.~Curtain, M.~Demetriou, and K.~Ito.
\newblock Adaptive observers for slowly time varying infinite dimensional
  systems.
\newblock In {\em Proceedings of the 37th IEEE Conference on Decision and
  Control (Cat. No.98CH36171)}, volume~4, pages 4022--4027 vol.4, 1998.

\bibitem{DEMETRIOU2018220}
M.~A. Demetriou.
\newblock Design of adaptive output feedback synchronizing controllers for
  networked pdes with boundary and in-domain structured perturbations and
  disturbances.
\newblock {\em Automatica}, 90:220--229, 2018.

\bibitem{DEMETRIOU19965346}
M.~A. Demetriou and K.~Ito.
\newblock Adaptive observers for a class of infinite dimensional systems.
\newblock {\em IFAC Proceedings Volumes}, 29(1):5346--5350, 1996.
\newblock 13th World Congress of IFAC, 1996, San Francisco USA, 30 June - 5
  July.

\bibitem{DECHcdc17}
G.~I. Detorakis and A.~Chaillet.
\newblock Incremental stability of spatiotemporal delayed dynamics and
  application to neural fields.
\newblock In {\em 56th IEEE Conference on Decision and Control}, pages
  5937--5942, 2017.

\bibitem{FARZA20092292}
M.~Farza, M.~M’Saad, T.~Maatoug, and M.~Kamoun.
\newblock Adaptive observers for nonlinearly parameterized class of nonlinear
  systems.
\newblock {\em Automatica}, 45(10):2292--2299, 2009.

\bibitem{Faugeras:2008wx}
O.~Faugeras, F.~Grimbert, and J.-J. Slotine.
\newblock {Absolute stability and complete synchronization in a class of neural
  fields models}.
\newblock {\em SIAM Journal of Applied Mathematics}, 61(1):205--250, 2008.

\bibitem{miyadera1992nonlinear}
I.~Miyadera.
\newblock {\em Nonlinear semigroups}, volume 109.
\newblock American Mathematical Soc., 1992.

\bibitem{potthast2009inverse}
R.~Potthast and P.~B. Graben.
\newblock Inverse problems in neural field theory.
\newblock {\em SIAM Journal on Applied Dynamical Systems}, 8(4):1405--1433,
  2009.

\bibitem{https://doi.org/10.48550/arxiv.2112.05497}
A.~Pyrkin, A.~Bobtsov, R.~Ortega, and A.~Isidori.
\newblock An adaptive observer for uncertain linear time-varying systems with
  unknown additive perturbations, 2021.

\bibitem{sastry1990adaptive}
S.~Sastry, M.~Bodson, and J.~F. Bartram.
\newblock Adaptive control: stability, convergence, and robustness, 1990.

\bibitem{showalter_monotone_2013}
R.~E. Showalter.
\newblock {\em Monotone operators in {Banach} space and nonlinear partial
  differential equations}, volume~49.
\newblock American Mathematical Society, 2013.

\bibitem{temam1986infinite}
R.~Temam.
\newblock Infinite-dimensional dynamical systems.
\newblock {\em Nonlinear functional analysis and its applications, Part 2},
  45(Part 2):431, 1986.

\bibitem{8786148}
J.~Wang, D.~Efimov, and A.~A. Bobtsov.
\newblock On robust parameter estimation in finite-time without persistence of
  excitation.
\newblock {\em IEEE Transactions on Automatic Control}, 65(4):1731--1738, 2020.

\bibitem{wilson1973mathematical}
H.~R. Wilson and J.~D. Cowan.
\newblock A mathematical theory of the functional dynamics of cortical and
  thalamic nervous tissue.
\newblock {\em Kybernetik}, 13(2):55--80, 1973.

\end{thebibliography}

\end{document}